\documentclass[12pt]{article}
\usepackage{textcomp}
\usepackage{graphicx}           
\usepackage{epstopdf}
\usepackage{color}              
\usepackage{xcolor,graphicx}
\usepackage{geometry}
\usepackage{float}
\usepackage{times}
\usepackage{indentfirst}        
\usepackage{amsmath,amssymb,bm} 
\usepackage{amsthm}
\usepackage{cases}              
\usepackage{setspace}
\usepackage{hyperref}
\usepackage{titlesec}
\usepackage[all]{xy}            
\usepackage{tikz}
\usepackage{multirow}
\usepackage{graphicx} 
\usepackage{epstopdf}
\usepackage{caption}
\usepackage{diagbox}
\usepackage{mathrsfs}
\usepackage{tikz}
\usepackage{dsfont}
\usepackage{graphics}
\usepackage{hyperref,cleveref,cite}
\usepackage{tcolorbox}
\usepackage{fancyhdr,fancyvrb}
\usepackage{xcolor,graphicx,geometry,float,indentfirst,cases,setspace,titlesec}
\usepackage{diagbox}

\newtheorem{thm}{Theorem}[section]

\newtheorem{lem}[thm]{Lemma}

\newtheorem{pf}{proof}[section]
\theoremstyle{remark}

\begingroup
\makeatletter \@addtoreset{equation}{section} \makeatother
\makeindex 
\def\qed{\hfill \rule{4pt}{7pt}}
\def\pf{\vskip 0.2cm {\noindent \bf Proof.}\quad}

\begin{document}

\begin{center}

 {\Large \bf Convexity and log-concavity of the partition function\\[5pt] weighted by the parity of the crank}

\end{center}

\begin{center}
Janet J.W. Dong$^{1}$ and {Kathy Q. Ji}$^{2}$ \vskip 2mm

$^{1,\,2}$ Center for Applied Mathematics\\[3pt]
Tianjin University\\[3pt]
Tianjin 300072, P.R. China\\[6pt]
   \vskip 2mm

 Emails: $^1$dongjinwei@tju.edu.cn  and $^2$kathyji@tju.edu.cn
\end{center}

\vskip 6mm \noindent {\bf Abstract.} Let $M_0(n)$ (resp. $M_1(n)$) denote the number of partitions of $n$ with even  (reps. odd) crank. Choi, Kang and Lovejoy
established   an asymptotic formula for $M_0(n)-M_1(n)$. By utilizing this formula with the explicit bound, we show that  $M_k(n-1)+M_k(n+1)>2M_k(n)$ for $k=0$ or $1$ and $n\geq 39$.  This result  can be seen as the refinement of the classical result regarding the convexity of  the partition function $p(n)$, which counts the number of partitions of $n$. We also show that $M_0(n)$ (resp. $M_1(n)$)  is
log-concave for $n\geq 94$ and satisfies the higher order Tur\'an inequalities for  $n\geq 207$ with the aid of  the upper bound and the lower bound for $M_0(n)$ and $M_1(n)$.



 \vskip 6mm

\section{Introduction}

This paper carries out a study of the partitions with even (resp. odd) crank from an analytic view.
 The crank of a partition
was defined by Andrews and Garvan \cite{Andrews-Garvan-1988} as the largest part if the partition contains no
ones, and otherwise as the number of parts larger than the number of ones minus the
number of ones. Let $p(n)$ denote the number of partitions of $n$. It is known that  the crank
provides combinatorial explanations for Ramanujan's
famous congruences $p(5n+4)\equiv 0 \pmod{5}$, $p(7n+5)\equiv 0 \pmod{7}$ and $p(11n+6)\equiv 0 \pmod{11}$, see Andrews and Garvan \cite{Andrews-Garvan-1988},  Dyson \cite{Dyson-1944} and Garvan \cite{Garvan-1988}. More precisely,  Let $M(r, Q ; n)$ be the number of partitions of $n$ with crank congruent to $r$ modulo $Q$. Andrews and Garvan \cite{Andrews-Garvan-1988}, building on the work of Garvan \cite{Garvan-1988}, established the following results:
\begin{align*}
M(r, 5 ; 5 n+4)&=\frac{p(5 n+4)}{5} \quad \text{for each}\quad 0\leq r\leq 4,\\[5pt]
M(r, 7 ; 7 n+5)&=\frac{p(7 n+5)}{7}\quad\text{for each} \quad  0\leq r\leq 6,\\[5pt]
M(r, 11 ; 11n+6)&=\frac{p(11 n+6)}{11}\quad\text{for each} \quad 0\leq r\leq 10.
\end{align*}
Recently,
Hamakiotes, Kriegman and Tsai \cite{Hamakiotes-Kriegman-Tsai-2021} demonstrated that $M(r, Q; n)$ is asymptotically equidistributed modulo odd number $Q$. More precisely, let $0 \leq r<Q$ with $Q$ an odd integer,    they established the following asymptotic result:
\[
\lim_{n\rightarrow \infty}\frac{M(r, Q ; n)}{p(n)} = \frac{1}{Q}.\]
Their proof relies on the asymptotic formula for $M(r, Q; n)$  derived by Zapata  Rol\'on \cite{Zapata-Rolon-2013}.

The study of partitions with even and odd cranks  was initially undertaken by Andrews and Lewis \cite{Andrews-Lewis-2000}.  For convenient, let $M_0(n)$ and $M_1(n)$ denote  the number of partitions of $n$ with even and odd cranks, respectively. Andrews and Lewis \cite{Andrews-Lewis-2000} showed that
\begin{thm}[Andrews-Lewis] \label{Andrews-Lewis} For $n\geq 0$,
\[(-1)^n (M_0(n)-M_1(n))>0.
\]
\end{thm}

Following Andrews and Lewis' footsteps, Choi, Kang and Lovejoy \cite{Choi-Kang-Lovejoy-2009} conducted a comprehensive study of $M_0(n)$ and $M_1(n)$. Among their main results, they derived  a family of Ramanujan type congruences satisfied by $M_0(n)-M_1(n)$ and obtained the following asymptotic formula for $M_0(n)-M_1(n)$.

\begin{thm}[Choi-Kang-Lovejoy]\label{choi-beta-thm}
	\begin{equation}\label{choi-formula-beta}
	M_0(n)-M_1(n)=E(n)+\frac{   \sqrt{6}\pi}{3 \mu(n)}\sum_{0<  j< \frac{\sqrt{3}\mu(n)}{2\sqrt{\pi}}}\cosh\left(\frac{\mu(n)}{2j}\right)\frac{\hat{A}_{j}(n)}{\sqrt{j}},
	\end{equation}
	where
	\begin{equation}\label{chern-beta-E}
		|E(n)|\leq \frac{95\cdot 6^{1/4}}{\sqrt{2\pi}}\mu(n)^{\frac{1}{2}}
	\end{equation}
	and
	\begin{equation}\label{defi-AM}
		\hat{A}_{j}(n)=\sum_{0\leq h< 2j\atop {\rm gcd}(h,2j)=1}\exp\left(-\frac{\pi i nh }{j}-\pi i\left(3s(h,2j)-2s\left(h,j\right)\right)\right),
	\end{equation}
	where  $s(h,j)$ is the  Dedekind sum   defined by
\begin{align}\label{def-DS}
s(h,j)= \left\{
	\begin{aligned}
			\sum_{r=1}^{j-1}\left(\frac{r}{j}-\left[\frac{r}{j}\right]-\frac{1}{2}\right)\left(\frac{h r}{j}-\left[\frac{h r}{j}\right]-\frac{1}{2}\right) \quad &\text{for} \quad j\geq 2,\\[5pt]
			0 \qquad\qquad\qquad \qquad\quad &\text{for}\quad j=1.
	\end{aligned}\right.
\end{align}
\end{thm}
Throughout this paper, we adopt the following notation as used in \cite{Lehmer-1939}
\begin{equation}\label{defi-x}
	\mu(n)=\frac{\pi\sqrt{24n-1}}{6}.
\end{equation}
The main objective of this paper is to employ the asymptotic formula \eqref{choi-formula-beta} for $M_0(n)-M_1(n)$ to investigate the convexity and log-concavity of $M_0(n)$ and $M_1(n)$. We first show that  $M_0(n)$ and $M_1(n)$  are convex when $n\geq 39$.

\begin{thm}\label{thm-convex}
    For $k=0$ or $1$ and $n\geq 39$,  \[M_k(n-1)+M_k(n+1)>2M_k(n).\]
    \end{thm}
It should be noted that Theorem \ref{thm-convex} can be viewed as the refinements of the  classical result involving the convexity of $p(n)$, that is, $p(n-1)+p(n+1)\geq2p(n)$ for $n\geq 2$, see Gupta \cite{Gupta-1978}  or Honsberger \cite[pp. 237--239]{Honsberger-1991}.

Furthermore, we establish the log-concavity of $M_0(n)$ and $M_1(n)$.

\begin{thm}\label{thm-LC} For $k=0$ or $1$ and $n\geq 94$,
	\[M_k(n)^2\geq M_k(n-1)M_k(n+1).\]
	To wit,  the sequences $\{M_0(n)\}_{n\geq 94}$  and $\{M_1(n)\}_{n\geq 93}$ are log-concave or satisfy the Tur\'an inequalities.
\end{thm}

To prove Theorem \ref{thm-LC} , we establish the upper bound and the lower bound of $M_k(n-1)M_k(n+1)/M_k(n)^2$.

\begin{thm}\label{thm-bound-Yk}
Let
\begin{equation}\label{defi-Yk}
	Y_k(n):=\frac{M_k(n-1)M_k(n+1)}{M_k(n)^2}.
\end{equation}
For $k=0$ or $1$ and $\mu(n)\geq 115 $,
	\begin{equation}\label{thm-bound-Y<}
		Y_k(n)<1-\frac{\pi^4}{9\mu(n)^3}+\frac{4\pi^4}{9\mu(n)^4}-\frac{\pi^4}{3\mu(n)^5}+\frac{\frac{\pi^8}{81}+5}{\mu(n)^6}
	\end{equation}
	and
	\begin{equation}\label{thm-bound-Y>}
		Y_k(n)>1-\frac{\pi^4}{9\mu(n)^3}+\frac{4\pi^4}{9\mu(n)^4}-\frac{\pi^4}{3\mu(n)^5}-\frac{60}{\mu(n)^6}.
	\end{equation}
\end{thm}

These bounds further allow us to show that $M_0(n)$ and $M_1(n)$ satisfy the  higher order Tur\'an inequalities.

\begin{thm}\label{thm-highturan} For $k=0$ or $1$ and $n\geq 207$,
	\begin{align}
		&4(M_{k}(n)^2-M_{k}(n-1)M_{k}(n+1))(M_{k}(n+1)^2-M_{k}(n)M_k(n+2))\nonumber\\[5pt]
		&\quad\geq (M_{k}({n})M_{k}(n+1)-M_{k}(n-1)M_{k}(n+2))^2.
	\end{align}
	That is, the sequences $\{M_0(n)\}_{n\geq207}$ and $\{M_1(n)\}_{n\geq206}$  satisfy the  higher order Tur\'an inequalities.
\end{thm}
It should be noted that the Tur\'an inequalities and the higher order Tur\'an inequalities for  the  partition function $p(n)$ and its variations have been extensively investigated recently, see, for example,  Bringmann, Kane, Rolen and Tripp \cite{Bringmann-Kane-Rolen-Tripp-2021}, Chen \cite{Chen-2017}, Chen, Jia and Wang \cite{Chen-Jia-Wang-2019}, DeSalvo and Pak \cite{DeSalvo-Pak-2015},  Griffin, Ono, Rolen and Zagier \cite{Griffin-Ono-Rolen-Zagier-2019} and  Ono, Pujahari and Rolen \cite{Ono-Pujahari-Rolen-2019}.
 In particular, Griffin, Ono, Rolen and Zagier \cite{Griffin-Ono-Rolen-Zagier-2019} showed that  $p(n)$ satisfies the order $d$  Tur\'an inequalities  for   sufficiently large $n$, confirming a conjectured of Chen, Jia and Wang \cite{Chen-Jia-Wang-2019}. For the definition of the order $d$  Tur\'an inequalities, please see Chen, Jia and Wang \cite{Chen-Jia-Wang-2019} or Griffin, Ono, Rolen and Zagier \cite{Griffin-Ono-Rolen-Zagier-2019}.

To conclude the introduction,   let us say a few words about  the order  $d$ Tur\'an inequalities for
$M_k(n)$.
In fact, based on \eqref{lim-Mk} and in view of  Theorem $3$ and Corollary $4$ in  \cite{Griffin-Ono-Rolen-Zagier-2019}, we could derive that

\begin{thm}
	When $k=0$ or $1$ and $d \geq 1$, $M_k(n)$ satisfies the order $d$ Tur\'an inequalities for sufficiently large $n$.
\end{thm}

It would be interesting to establish the minimal number $O_{M_k}(d)$ such that $M_k(n)$ satisfies the order $d$ Tur\'an inequalities for $n\geq O_{M_k}(d)$.  By Theorem \ref{thm-LC}, we see that
\[O_{M_0}(2)=93 \quad \text{and} \quad  O_{M_1}(2)=92.\]
From   Theorem \ref{thm-highturan}, we find that
\[O_{M_0}(3)=206 \quad \text{and} \quad  O_{M_1}(3)=205.\]

The paper is organized as follows. In Section 2, we first establish the error bound for the asymptotic formula for $M_0(n)-M_1(n)$ due to Choi, Kang and Lovejoy (that is, Theorem \ref{formula-beta-thm}).  We then give two additional application of Theorem \ref{formula-beta-thm}, see Theorem \ref{Andrews-Lewis-d} and Theorem \ref{coro-eq-Mp}. In Section 3, we employ   Theorem \ref{formula-beta-thm} to establish the upper bound and the lower bound for $M_0(n)$ and $M_1(n)$ (that is, Theorem \ref{thm-Mk-eq}), which is useful in the study of the convexity and the log-concavity  of $M_0(n)$ and $M_1(n)$. Section 4 focuses on exploring    the convexity of $M_0(n)$ and $M_1(n)$ by utilizing  Theorem \ref{thm-Mk-eq}.   Section 5 is devoted to employing Theorem \ref{thm-Mk-eq} to establish the upper bound and the lower bound of $M_k(n-1)M_k(n+1)/M_k(n)^2$ (that is, Theorem \ref{thm-bound-Yk}). These bounds play a crucial role  in the proofs that $M_0(n)$ and $M_1(n)$ are
log-concave for $n\geq 94$ and satisfy the higher order Tur\'an inequalities for  $n\geq 207$, which are established in  Section 6.

\section{The error bound}

 In this section, we state the error bound for the asymptotic formula of $M_0(n)-M_1(n)$ due to Choi, Kang and Lovejoy, which is required in our study of the convexity and log-concavity  of $M_0(n)$ and  $M_1(n)$. We also provide two additional applications of this error bound. First application is to  provide a direct analytic proof of Andrews and Lewis'  result (that is, Theorem \ref{Andrews-Lewis}) and the second one is to demonstrate that the cranks are asymptotically equidistributed modulo 2.

Using  Theorem \ref{choi-beta-thm}, we obtain the  following asymptotic formula for $M_0(n)-M_1(n)$ with an effective bound on the error term.

\begin{thm}\label{formula-beta-thm}
	For $\mu(n)\geq 4$, or equivalently, $n\geq 3$,
	\begin{equation}\label{formula-beta}
		M_0(n)-M_1(n)=\frac{(-1)^n\pi}{\sqrt{6}\mu(n)}e^{\frac{\mu(n)}{2}}+E_\beta(n),
	\end{equation}
	where
	\begin{equation}\label{bound-R}
		|E_\beta(n)|\leq 63 \mu(n)^{\frac{1}{2}}e^{\frac{\mu(n)}{4}}.
	\end{equation}
\end{thm}

 \pf Observing that $\hat{A}_1(n)=(-1)^n$, and using the formula of $\cosh(s)$ from \cite[p. 459]{Abra-1972}
\begin{equation}
	\cosh(s)=\frac{e^{s}+e^{-s}}{2},
\end{equation}
we deduce from
\eqref{choi-formula-beta} that for $\mu(n)\geq 4$,
\[M_0(n)-M_1(n)=\frac{(-1)^n\pi}{\sqrt{6}\mu(n)}e^{\frac{\mu(n)}{2}}+E_{\beta}(n),\]
where
\begin{equation}\label{eq-E-beta}
	E_{\beta}(n)=E(n)+\frac{(-1)^n\pi}{\sqrt{6}\mu(n)}e^{-\frac{\mu(n)}{2}}+\frac{   \sqrt{6}\pi}{3 \mu(n)}\sum_{2\leq j<\frac{\sqrt{3}\mu(n)}{2\sqrt{\pi}}}\cosh\left(\frac{\mu(n)}{2j}\right)\frac{\hat{A}_j(n)}{\sqrt{j}}.
\end{equation}
We next establish the  bound for $|E_{\beta}(n)|$. By the definition \eqref{defi-AM} of  $\hat{A}_j(n)$, we derive that  for any $n\geq 0$ and $j\geq 1$,
\begin{equation*}
	|\hat{A}_j(n)|\leq 2j.
\end{equation*}
Hence, we have
\begin{align}\label{eq-E-beta-1}
	\frac{   \sqrt{6}\pi}{3 \mu(n)}\sum_{2\leq  j<\frac{\sqrt{3}\mu(n)}{2\sqrt{\pi}}}\cosh\left(\frac{\mu(n)}{2j}\right)\frac{\hat{A}_j(n)}{\sqrt{j}}
 &\leq \frac{2\sqrt{6}\pi}{3 \mu(n)}\sum_{2\leq j\leq\frac{\mu(n)}{2}} \sqrt{j}\cosh \left(\frac{\mu(n)}{2j}\right)\nonumber\\[5pt]
 &\leq\frac{ 2\pi}{\sqrt{3} \mu(n)^{\frac{1}{2}}}\sum_{2\leq j\leq\frac{\mu(n)}{2}} \cosh\left(\frac{\mu(n)}{2j}\right)\nonumber\\[5pt]
 &\leq \frac{\pi \mu(n)^{\frac{1}{2}}}{\sqrt{3}} \cosh\left(\frac{\mu(n)}{4}\right)\nonumber\\[5pt]
 &\leq \frac{\pi \mu(n)^{\frac{1}{2}}}{2\sqrt{3}} \left(e^{\frac{\mu(n)}{4}}+e^{-\frac{\mu(n)}{4}}\right).
\end{align}
Applying \eqref{chern-beta-E} and \eqref{eq-E-beta-1} to  \eqref{eq-E-beta} , we derive that for $\mu(n)\geq 4$,
\begin{align}
	|E_{\beta}(n)|&\le\frac{\pi}{\sqrt{6}\mu(n)}e^{-\frac{\mu(n)}{2}}+\frac{\pi \mu(n)^{\frac{1}{2}}}{2\sqrt{3}} \left(e^{\frac{\mu(n)}{4}}+e^{-\frac{\mu(n)}{4}}\right)+\frac{95\cdot 6^{1/4}}{\sqrt{2\pi}}\mu(n)^{\frac{1}{2}}\nonumber\\[5pt]
	&\le \frac{  \pi}{2\sqrt{3} }\mu(n)^{\frac{1}{2}}e^{\frac{\mu(n)}{4}}+\left(\frac{   \pi}{4\sqrt{6}}+\frac{   \pi}{2\sqrt{3}}+60 \right)\mu(n)^{\frac{1}{2}}\nonumber\\[5pt]
	&\le \left(\frac{\pi}{4\sqrt{6}}+\frac{\pi}{\sqrt{3}}+60 \right)\mu(n)^{\frac{1}{2}} e^{\frac{\mu(n)}{4}} \nonumber\\[5pt]
	&\le 63 \mu(n)^{\frac{1}{2}} e^{\frac{\mu(n)}{4}}.\nonumber
\end{align}
This completes the proof.
\qed

Using Theorem \ref{formula-beta-thm}, we give a direct analytic proof of Theorem \ref{Andrews-Lewis}. In particular, we obtain the following more general result.

\begin{thm}\label{Andrews-Lewis-d} Given a positive integer $d\geq 1$. For $n\geq \left\lceil\frac{24}{\pi^2}\left(\ln\left(\frac{7d}{2}\right)\right)^2+\frac{1}{24}\right\rceil$,
\begin{equation}\label{Andrews-Lewis-d-eq}
(-1)^n (M_0(n)-M_1(n))> d.
\end{equation}
\end{thm}

\pf Using Theorem \ref{formula-beta-thm}, we find that for $n\geq 3$,
\begin{align} \label{formula-beta-thmaaa}
    (-1)^n(M_0(n)-M_1(n))=\frac{\pi}{\sqrt{6}\mu(n)}e^{\frac{\mu(n)}{2}}+(-1)^nE_\beta(n),
\end{align}
where
\begin{align*}
\left|E_\beta(n)\right|\leq  63 \mu(n)^{\frac{1}{2}}e^{\frac{\mu(n)}{4}}.
\end{align*}
We first show that for $n\geq 1$,
\begin{align}\label{eq-AL-2-1}
 (-1)^n(M_0(n)-M_1(n))>\frac{1}{7}\mu(n)^{\frac{1}{2}}e^{\frac{\mu(n)}{4}}.
\end{align}
By \eqref{formula-beta-thmaaa},  it  suffices to show that for $n\geq 1$,
\begin{align}\label{eq-AL-2}
\frac{\pi}{\sqrt{6}\mu(n)}e^{\frac{\mu(n)}{2}}- 63 \mu(n)^{\frac{1}{2}}e^{\frac{\mu(n)}{4}}>\frac{1}{7}\mu(n)^{\frac{1}{2}}e^{\frac{\mu(n)}{4}}.
\end{align}
We claim that for $\mu(n)\geq 38$, or equivalently, $n\geq 220$,
\begin{equation}\label{eq-AL-3}
    \frac{7\pi}{\sqrt{6}\mu(n)}e^{\frac{\mu(n)}{4}}>442 \mu(n)^{\frac{1}{2}}.
\end{equation}
Define
\[f(s)=\frac{7\pi}{442\sqrt{6}s^{\frac{3}{2}}}e^{\frac{s}{4}}.\]
It is easy to derive that for $s\geq 6$,
\begin{equation*}
    f'(s)=\frac{7\pi  (s-6)e^{\frac{s}{4}}}{1768 \sqrt{6} s^{\frac{5}{2}}}\geq 0.
\end{equation*}
Hence for $\mu(n)\geq 38$,
\begin{equation*}
    f(\mu(n))=\frac{7\pi}{442\sqrt{6}\mu(n)^{\frac{3}{2}}}e^{\frac{\mu(n)}{4}}\geq f(38)>1,
\end{equation*}
and  so \eqref{eq-AL-3}   is proved. Using  \eqref{eq-AL-3},  we derive that  \eqref{eq-AL-2} is valid for $n\geq 220$. It follows that  \eqref{eq-AL-2-1} is valid for $n\geq 220$. It can be checked that \eqref{eq-AL-2-1} is  valid for $1 \leq n\leq 220$. Hence \eqref{eq-AL-2-1} is proved. On the other hand, when $d\geq 1$ and for  $\mu(n)\geq 4\ln
\left(\frac{7d}{2}\right)$, or equivalently,  $n\geq \left\lceil\frac{24}{\pi^2}\left(\ln\left(\frac{7d}{2}\right)\right)^2+\frac{1}{24}\right\rceil$, we see that
\begin{equation}\label{eq-AL-4}
    \frac{1}{7}\mu(n)^{\frac{1}{2}}e^{\frac{\mu(n)}{4}}\geq\frac{1}{7}\sqrt{4\ln
\left(\frac{7d}{2}\right)}e^{\frac{4\ln
\left(\frac{7d}{2}\right)}{4}} >\frac{2}{7}\cdot \frac{7d}{2} =d.
\end{equation}
Combining \eqref{eq-AL-2-1} and \eqref{eq-AL-4}, we obtain  \eqref{Andrews-Lewis-d-eq}. This  completes the proof of Theorem \ref{Andrews-Lewis-d}.
\qed

We conclude this section by showing  that the cranks are asymptotically equidistributed modulo $2$.
\begin{thm}\label{coro-eq-Mp}
	For $k=0$ or $1$,
\begin{equation}\label{lim-Mk}
	\lim_{n\rightarrow \infty}\frac{M_k(n)}{p(n)}=\frac{1}{2}.
\end{equation}
\end{thm}
It should be noted that asymptotically  equidistribution of partition ranks modulo $2$ has been recently proved by   Gomez and Zhu \cite{Gomez-Zhu-2021} and Masri \cite{Masri-2016, Masri-2021}.

 It turns out   the proof of   Theorem \ref{coro-eq-Mp} also requires the following  lower bound for $p(n)$ given by   Bessenrodt and Ono \cite{Bessenrodt-Ono-2016}:  For $n \geq 1$,
\begin{equation} \label{BO}
p(n)>\frac{\sqrt{3}}{12 n}\left(1-\frac{1}{\sqrt{n}}\right) e^{\mu(n)}.
\end{equation}

We aim to prove the following consequence,  which leads to the immediate proof of  Theorem \ref{coro-eq-Mp}.

\begin{thm}\label{thm-eq-Mp}
	For $k=0$ or $1$ and for $n\geq 4$,
	\begin{equation}\label{thm-eq-Mpmt}
	\frac{M_k(n)}{p(n)}=\frac{1}{2}+(-1)^kE^{c}(n),
	\end{equation}
	where
	\begin{equation}\label{thm-eq-Mperror}
	\left|E^{c}(n)\right|\leq 11578 e^{-\frac{\mu(n)}{4}}.
	\end{equation}
\end{thm}

\pf  By definition, we see that
\begin{equation}\label{mp2}
	p(n)=M_0(n)+M_1(n).
\end{equation}
It follows that for $k=0$ or $1$,
\[
\frac{M_k(n)}{p(n)}=\frac{1}{2}+(-1)^k \frac{M_0(n)-M_1(n)}{2 p(n)}. \]
Assume that
\[E^c(n):=\frac{M_0(n)-M_1(n)}{2 p(n)}.\]
In light of Theorem  \ref{formula-beta-thm} and \eqref{BO},   we derive that for $n\geq 4$ ,
\begin{align}
	\left|E^c(n)\right| & \leq \frac{12\mu(n)^2}{\sqrt{3}} e^{-\mu(n)}\left(\frac{\pi}{\sqrt{6}\mu(n)}e^{\frac{\mu(n)}{2}}+|E_{\beta}(n)|\right) \nonumber \\[5pt]
	& \leq 2\sqrt{2}\pi \mu(n) e^{-\frac{\mu(n)}{2}}+\frac{756}{\sqrt{3}}\mu(n)^{\frac{5}{2}} e^{-\frac{3\mu(n)}{4}} \nonumber \\[5pt]
	& \leq\left(2\sqrt{2}\pi  +\frac{756}{\sqrt{3}}\right) \mu(n)^{\frac{5}{2}} e^{-\frac{\mu(n)}{2}} \nonumber \\[5pt]
	& \leq446\mu(n)^{\frac{5}{2}} e^{-\frac{\mu(n)}{2}}. \nonumber
\end{align}
We claim that for  $\mu(n)>0$,
\begin{equation}\label{cor-claim}
	446\mu(n)^{\frac{5}{2}} e^{-\frac{\mu(n)}{4}}<11578.
\end{equation}
Define
\[m(s):=446s^{\frac{5}{2}} e^{-\frac{s}{4}}.\]
It is evident that
\[m'(s)=\frac{223}{2}s^{\frac{3}{2}} (10-s)e^{-\frac{s}{4}}.\]
Since $m'(s)<0$ for $s>10$ and $m'(s)>0$ for $0<s<10$, we derive that  $m(s)$  attains its maximum value     at $s=10$, so
\[m(\mu(n))\leq m(10)<11578,\]
and hence \eqref{cor-claim} holds. We therefore obtain \eqref{thm-eq-Mperror}.  This completes the proof. \qed

\section{An upper bound and a lower bound for $M_k(n)$}

In this section, we aim to establish  the following upper bound and the lower bound for $M_0(n)$ and $M_1(n)$ (Theorem \ref{thm-Mk-eq}). It turns out that the proof of Theorem \ref{thm-Mk-eq} also requires the following  effective bound on $p(n)$  due to Locus Dawsey and Masri  \cite[Lemma 4.2]{LocusDawsey-Masri-2019}.  For $n\geq 1$,
	\begin{align} \label{eq-p}
		p(n)=\frac{\pi^2}{6\sqrt{3}\mu(n)^2}\left(1-\frac{1}{\mu(n)}\right)e^{\mu(n)}+E_p(n),
  \end{align}
	where
	\begin{equation}\label{bound-RR}
	\left|E_p(n)\right| \leq1313e^{\frac{\mu(n)}{2}}.
	\end{equation}

\begin{thm}\label{thm-Mk-eq}
	Let \begin{equation}\label{defi-M}
		G(n):=\frac{\pi^2}{12\sqrt{3}\mu(n)^2}\left(1-\frac{1}{\mu(n)}\right)e^{\mu(n)}.
	\end{equation}
	Then for $k=0,1$ and $\mu(n)\geq88$,
	\begin{align}\label{eq-Mk}
		G(n)\left(1-\frac{1}{\mu(n)^6}\right)\le M_k(n)\le G(n)\left(1+\frac{1}{\mu(n)^6}\right).
	\end{align}
\end{thm}

\pf  We first show that when $k=0$ or $1$ and for   $n\geq 3$,
	\begin{align} \label{thm-Mk-eq-pfaa}
		M_k(n)=\frac{\pi^2}{12\sqrt{3}\mu(n)^2}\left(1-\frac{1}{\mu(n)}\right)e^{\mu(n)}+R^c_{k}(n),	
	\end{align}
	where
	\begin{equation}
		|R^c_{k}(n)|\le 689 e^{\frac{\mu(n)}{2}}.
	\end{equation}
  Applying \eqref{mp2} to  \eqref{eq-p} and using Theorem \ref{formula-beta-thm},
 we could derive  that for $k=0$ or $1$ and  $\mu(n)\geq 4$,
\begin{align}
	M_k(n) =\frac{\pi^2}{12\sqrt{3}\mu(n)^2}\left(1-\frac{1}{\mu(n)}\right)e^{\mu(n)}+R_k^c(n),\nonumber
\end{align}
where
\[
R_k^c(n)=\frac{(-1)^{n+k}\pi}{2\sqrt{6}\mu(n)}e^{\frac{\mu(n)}{2}}+\frac{1}{2}\left(E_p(n)+(-1)^kE_{\beta}(n)\right).
\]
In light of \eqref{bound-R} and \eqref{bound-RR}, we see that   for $k=0$ or $1$ and $\mu(n)\geq 4$,
\begin{align*}
	\left|R_k^c(n)\right|
	&\le\frac{\pi}{2\sqrt{6}\mu(n)}e^{\frac{\mu(n)}{2}}+\frac{1313}{2}e^{\frac{\mu(n)}{2}}+\frac{63}{2} \mu(n)^{\frac{1}{2}}e^{\frac{\mu(n)}{4}}.	 \\[5pt]
	&\le \left(\frac{\pi}{8\sqrt{6}}+\frac{1313}{2}+\frac{63}{2}\right)e^{\frac{\mu(n)}{2}} \\[5pt]
	&\le 689e^{\frac{\mu(n)}{2}}.
\end{align*}	
This completes the proof of  \eqref{thm-Mk-eq-pfaa}.

Define
\begin{align}
	T(n):=\frac{689 e^{\frac{\mu(n)}{2} }}{\frac{\pi^2}{12\sqrt{3}\mu(n)^2}\left(1-\frac{1}{\mu(n)}\right)e^{\mu(n)}}
	=\frac{8268\sqrt{3}\mu(n)^2}{\pi^2\left(1-\frac{1}{\mu(n)}\right)}e^{-\frac{\mu(n)}{2}}.
\end{align}
Using  \eqref{thm-Mk-eq-pfaa}, we find that for $n\geq 3$,
\[G(n)\left(1-T(n)\right)\le M_k(n)\le G(n)\left(1+T(n)\right).\]
To show \eqref{eq-Mk}, it is enough to prove that for $\mu(n)\geq 88$,
\begin{align}\label{eq-G1}
	T(n)\le \frac{1}{\mu(n)^6}.
\end{align}
Note that for $\mu(n)\ge 2$,
\[\left(1-\frac{1}{\mu(n)}\right)\left(1+\frac{2}{\mu(n)}\right)=1+\mu(n)^{-2}(\mu(n)-2)\ge 1.\]
Hence for $n\geq 3$,
\begin{align}\label{eq-G2}
	T(n)\le\frac{8268\sqrt{3}}{\pi^2}\left(\mu(n)^2+2\mu(n)\right)e^{-\frac{\mu(n)}{2}}.
\end{align}
We claim that for $\mu(n)\geq 88 $,
\begin{align}\label{eq-G3}
	\frac{8268\sqrt{3}}{\pi^2}e^{-\frac{\mu(n)}{2}}\le\frac{1}{2\mu(n)^8}.
\end{align}
Define
\[L(s):=\frac{16536\sqrt{3}}{\pi^2}s^8e^{-\frac{s}{2}}.\]
It is evident that
\[L'(s)=\frac{8268\sqrt{3}}{\pi^2} e^{-\frac{s}{2}} (16-s) s^7.\]
Since $L'(s)\le 0$ for $s\ge 16$, we deduce that $L(s)$ is decreasing when $s\geq 16$, this implies that
\[L(\mu(n))=\frac{16536\sqrt{3}}{\pi^2}\mu(n)^8e^{-\frac{\mu(n)}{2}}\le L(88)<1 \]
for $\mu(n)\geq 88$, and so the claim is proved. Applying \eqref{eq-G3} to \eqref{eq-G2}, we obtain \eqref{eq-G1}. This completes the proof.\qed

\section{The convexity of $M_0(n)$ and $M_1(n)$}
The main objective of this section is   to establish the convexity of $M_0(n)$ and $M_1(n)$ (Theorem \ref{thm-convex}). We begin by  proving the following two inequalities.

\begin{lem}\label{lem-convex}
For $\mu(n)\geq 6$,
\begin{align}\label{lem-cx-1}
    \frac{G(n-1)}{G(n)}>\left(1+\frac{2\pi^2}{3\mu(n)^2}\right)\left(1-\frac{\pi^2}{\mu(n)^3}\right)\left(1-\frac{\pi^2}{3\mu(n)}+\frac{\pi^4}{18\mu(n)^2}-\frac{18\pi^4+\pi^6}{162\mu(n)^3}\right)
\end{align}
and
\begin{align}\label{lem-cx-2}
    \frac{G(n+1)}{G(n)}>\left(1-\frac{2\pi^2}{3\mu(n)^2}\right)\left(1+\frac{\pi^2}{3\mu(n)^3}\right)\left(1+\frac{\pi^2}{3\mu(n)}+\frac{\pi^4}{18\mu(n)^2}-\frac{\pi^4}{9\mu(n)^3}\right).
\end{align}
\end{lem}
\proof
Recall that \[
		G(n)=\frac{\pi^2}{12\sqrt{3}\mu(n)^2}\left(1-\frac{1}{\mu(n)}\right)e^{\mu(n)}.\]
We have
\begin{align}\label{lem-cx-3}
    \frac{G(n-1)}{G(n)}=\frac{\mu(n)^2}{\mu(n-1)^2}\cdot \frac{1-\frac{1}{\mu(n-1)}}{1-\frac{1}{\mu(n)}}\cdot \exp\left(\mu(n-1)-\mu(n)\right)
\end{align}
and
\begin{align}\label{lem-cx-4}
    \frac{G(n+1)}{G(n)}=\frac{\mu(n)^2}{\mu(n+1)^2}\cdot \frac{1-\frac{1}{\mu(n+1)}}{1-\frac{1}{\mu(n)}}\cdot \exp\left(\mu(n+1)-\mu(n)\right).
\end{align}
From the definition \eqref{defi-x} of $\mu(n)$, we see that
for $\mu(n)\geq 3$,
\begin{equation}\label{lem-cx-aaaa}
\mu(n-1)^2=\mu(n)^2-\frac{2\pi^2}{3} \quad \text{and} \quad \mu(n+1)^2=\mu(n)^2+\frac{2\pi^2}{3}. \end{equation}
It is evident that
\begin{equation*}
 \mu(n)^2-\left(\mu(n)^2-\frac{2\pi^2}{3} \right)\left(1+\frac{2\pi^2}{3\mu(n)^2} \right)=\frac{4\pi^4}{9\mu(n)^2}>0
\end{equation*}
and
\begin{equation*}
 \mu(n)^2-\left(\mu(n)^2+\frac{2\pi^2}{3} \right)\left(1-\frac{2\pi^2}{3\mu(n)^2} \right)=\frac{4\pi^4}{9\mu(n)^2}>0.
\end{equation*}
Hence, we derive that for $\mu(n)\geq 3$,
\begin{equation}\label{lem-cx-5}
  \frac{\mu(n)^2}{\mu(n-1)^2}>1+\frac{2\pi^2}{3\mu(n)^2}
\end{equation}
and
\begin{equation}\label{lem-cx-6}
  \frac{\mu(n)^2}{\mu(n+1)^2}>1-\frac{2\pi^2}{3\mu(n)^2},
\end{equation}

We proceed to show that for $\mu(n)\geq 6$,
\begin{align}\label{lem-cx-7}
    \frac{1-\frac{1}{\mu(n-1)}}{1-\frac{1}{\mu(n)}}>1-\frac{\pi^2}{\mu(n)^3}
\end{align}
and
\begin{align}\label{lem-cx-8}
   \frac{1-\frac{1}{\mu(n+1)}}{1-\frac{1}{\mu(n)}}>1+\frac{\pi^2}{3\mu(n)^3}.
\end{align}
By \eqref{lem-cx-aaaa}, we see that for $\mu(n)\geq 3$,
\begin{equation}\label{lem-cx-aaaacc}
\mu(n-1)=\sqrt{\mu(n)^2-\frac{2 \pi^2}{3}} \quad \text{and} \quad   \mu(n+1)=\sqrt{\mu(n)^2+\frac{2 \pi^2}{3}} .
\end{equation}
Then we have
\begin{align} \label{lem-cx-aaaaccdd}
\begin{aligned}
& \mu(n-1)=\mu(n)-\frac{\pi^2}{3 \mu(n)}-\frac{\pi^4}{18 \mu(n)^3}-\frac{\pi^6}{54 \mu(n)^5}-\frac{5 \pi^8}{648 \mu(n)^7}+o\left(\frac{1}{\mu(n)^8}\right), \\[5pt]
& \mu(n+1)=\mu(n)+\frac{\pi^2}{3 \mu(n)}-\frac{\pi^4}{18 \mu(n)^3}+\frac{\pi^6}{54 \mu(n)^5}-\frac{5 \pi^8}{648 \mu(n)^7}+o\left(\frac{1}{\mu(n)^8}\right).
\end{aligned}
\end{align}
It can be checked that for $\mu(n)\geq 6$,
\begin{equation}\label{lem-cx-9}
    \mu(n-1)>w(n)\quad \text{and} \quad \mu(n+1)>y(n),
\end{equation}
where
\[w(n)=\mu(n)-\frac{\pi^2}{3\mu(n)}-\frac{\pi^4}{9\mu(n)^3}\]
and
\[y(n)=\mu(n)+\frac{\pi^2}{3\mu(n)}-\frac{\pi^4}{18\mu(n)^3}.\]
Hence, we obtain that for $\mu(n)\geq 6$,
\begin{align}\label{lem-cx-10}
    \frac{1-\frac{1}{\mu(n-1)}}{1-\frac{1}{\mu(n)}}>\frac{1-\frac{1}{w(n)}}{1-\frac{1}{\mu(n)}}
\end{align}
and
\begin{align}\label{lem-cx-11}
   \frac{1-\frac{1}{\mu(n+1)}}{1-\frac{1}{\mu(n)}}>\frac{1-\frac{1}{y(n)}}{1-\frac{1}{\mu(n)}}.
\end{align}
To prove \eqref{lem-cx-7} and \eqref{lem-cx-8}, it is enough to show that for $\mu(n)\geq 6$,
\begin{align}\label{lem-cx-12}
    \frac{1-\frac{1}{w(n)}}{1-\frac{1}{\mu(n)}}-\left(1-\frac{\pi^2}{\mu(n)^3}\right)>0
\end{align}
and
\begin{align}\label{lem-cx-13}
   \frac{1-\frac{1}{y(n)}}{1-\frac{1}{\mu(n)}}-\left(1+\frac{\pi^2}{3\mu(n)^3}\right)>0.
\end{align}
Observe that
\begin{align}\label{lem-cx-14}
    \frac{1-\frac{1}{w(n)}}{1-\frac{1}{\mu(n)}}-\left(1-\frac{\pi^2}{\mu(n)^3}\right)=\frac{\phi_1(\mu(n))}{\mu(n)^6(\mu(n)-1)w(n)},
\end{align}
where
\[\phi_1(s)=\frac{2\pi^2}{3}s^5-\pi^2s^4-\frac{4\pi^4}{9}s^3+\frac{\pi^4}{3}s^2-\frac{\pi^6}{9}s+\frac{\pi^6}{9}.\]
It can be checked that for $s\geq 4$,
\begin{align*}\left\{
	\begin{aligned}
		&	\frac{2\pi^2}{3}s^5-\pi^2s^4-\frac{4\pi^4}{9}s^3>0,\\[5pt]
		&	\frac{\pi^4}{3}s^2-\frac{\pi^6}{9}s+\frac{\pi^6}{9}>0,
	\end{aligned}\right.
\end{align*}
which implies that $\phi_1(s)>0$ for $s\geq 4$, and so $ \phi_1(\mu(n))>0$ for  $\mu(n)\geq 4$. Hence \eqref{lem-cx-12} holds.

Similarly, we can write
\begin{align}\label{lem-cx-16}
    \frac{1-\frac{1}{y(n)}}{1-\frac{1}{\mu(n)}}-\left(1+\frac{\pi^2}{3\mu(n)^3}\right)=\frac{\phi_2(\mu(n))}{\mu(n)^6(\mu(n)-1)y(n)},
\end{align}
where
\[\phi_2(s)=\frac{\pi^2}{3}s^4-\frac{\pi^4}{6}s^3+\frac{\pi^4}{9}s^2+\frac{\pi^6}{54}s-\frac{\pi^6}{54}.\]
It can be checked that for $s\geq 5$,
\[\frac{\pi^2}{3}s^4-\frac{\pi^4}{6}s^3-\frac{\pi^6}{54}>0,\]
which implies that $\phi_2(s)>0$ for $s\geq 5$. It follows that $\phi_2(\mu(n))>0$ for $\mu(n)\geq 5$.
Hence \eqref{lem-cx-13} is valid.

We proceed to estimate $\exp\left(\mu(n-1)-\mu(n)\right)$ and $\exp\left(\mu(n+1)-\mu(n)\right)$. We claim that for $\mu(n)\geq 6$,
\begin{equation}\label{lem-cx-18}
    \exp\left(\mu(n-1)-\mu(n)\right)>1-\frac{\pi^2}{3\mu(n)}+\frac{\pi^4}{18\mu(n)^2}-\frac{18\pi^4+\pi^6}{162\mu(n)^3}.
\end{equation}
and
\begin{equation}\label{lem-cx-19}
    \exp\left(\mu(n+1)-\mu(n)\right)>1+\frac{\pi^2}{3\mu(n)}+\frac{\pi^4}{18\mu(n)^2}-\frac{\pi^4}{9\mu(n)^3}.
\end{equation}
Using \eqref{lem-cx-9}, we derive that for $\mu(n)\geq 6$,
\begin{equation}\label{lem-cx-20}
    \exp\left(\mu(n-1)-\mu(n)\right)>\exp\left(-\frac{\pi^2}{3\mu(n)}-\frac{\pi^4}{9\mu(n)^3}\right)
\end{equation}
and
\begin{equation}\label{lem-cx-21}
    \exp\left(\mu(n+1)-\mu(n)\right)>\exp\left(\frac{\pi^2}{3\mu(n)}-\frac{\pi^4}{18\mu(n)^3}\right).
\end{equation}
Observe that for $s<0$,
\begin{equation}\label{lem-cx-22}
    e^s>1+s+\frac{s^2}{2}+\frac{s^3}{6}
\end{equation}
and for $s>0$,
\begin{equation}\label{lem-cx-23}
    e^s>1+s+\frac{s^2}{2}.
\end{equation}
Hence, by \eqref{lem-cx-22}, we derive that for $\mu(n)\geq 4$,
\begin{align}
    \exp\left(-\frac{\pi^2}{3\mu(n)}-\frac{\pi^4}{9\mu(n)^3}\right)&>1-\frac{\pi^2}{3\mu(n)}+\frac{\pi^4}{18\mu(n)^2}-\frac{18\pi^4+\pi^6}{162\mu(n)^3}+\frac{\pi^6}{27\mu(n)^4}\nonumber\\[5pt]
    &\quad-\frac{\pi^8}{162\mu(n)^5}+\frac{\pi^8}{162\mu(n)^6}-\frac{\pi^{10}}{486\mu(n)^7}-\frac{\pi^{12}}{4374\mu(n)^9}.\nonumber\\[5pt]
    &>1-\frac{\pi^2}{3\mu(n)}+\frac{\pi^4}{18\mu(n)^2}-\frac{18\pi^4+\pi^6}{162\mu(n)^3}. \label{lem-cx-24}
\end{align}
The second inequality follows from the following observation: For $\mu(n)\geq 4$,
\begin{align*}\left\{
	\begin{aligned}
		&	\frac{\pi^6}{27\mu(n)^4}-\frac{\pi^8}{162\mu(n)^5}> 0,\\[5pt]
		&	\frac{\pi^8}{162\mu(n)^6}-\frac{\pi^{10}}{486\mu(n)^7}-\frac{\pi^{12}}{4374\mu(n)^9}> 0.
	\end{aligned}\right.
\end{align*}
Applying \eqref{lem-cx-24} to  \eqref{lem-cx-20}, we obtain \eqref{lem-cx-18}.

Using \eqref{lem-cx-23}, we derive that for $\mu(n)\geq 4$,
\begin{align}
    \exp\left(\frac{\pi^2}{3\mu(n)}-\frac{\pi^4}{18\mu(n)^3}\right)&>1+\frac{\pi^2}{3\mu(n)}+\frac{\pi^4}{18\mu(n)^2}-\frac{\pi^4}{18\mu(n)^3}-\frac{\pi^6}{54\mu(n)^4}+\frac{\pi^8}{648\mu(n)^6}.\nonumber\\[5pt]
    &>1+\frac{\pi^2}{3\mu(n)}+\frac{\pi^4}{18\mu(n)^2}-\frac{\pi^4}{9\mu(n)^3}. \label{lem-cx-25}
\end{align}
The second inequality is derived from the following observation, that is,  for $\mu(n)\geq 4$,
\[\frac{\pi^4}{18\mu(n)^3}-\frac{\pi^6}{54\mu(n)^4}>0.\]
 Combining \eqref{lem-cx-21} and \eqref{lem-cx-25} yields
\eqref{lem-cx-19}.

Applying \eqref{lem-cx-5}, \eqref{lem-cx-7} and \eqref{lem-cx-18} to \eqref{lem-cx-3}, we obtain  \eqref{lem-cx-1}. Substituting \eqref{lem-cx-6}, \eqref{lem-cx-8} and \eqref{lem-cx-19} into \eqref{lem-cx-4}, we obtain \eqref{lem-cx-2}. Thus we complete the proof of Lemma \ref{lem-convex}. \qed

We are now in a position to give a proof of Theorem \ref{thm-convex}.

\noindent{\it Proof of Theorem \ref{thm-convex}. }
We aim to prove that
\begin{equation}\label{convex-eq}
    \frac{M_k(n-1)+M_k(n+1)}{M_k(n)}>2.
\end{equation}
In light of Theorem \ref{thm-Mk-eq}, we deduce that for $\mu(n)\geq 88$,
\begin{align} \label{convex-tt}
   \frac{M_k(n-1)+M_k(n+1)}{M_k(n)} &\geq \frac{G(n-1)\left(1-\frac{1}{\mu(n-1)^4}\right)+G(n+1)\left(1-\frac{1}{\mu(n+1)^4}\right)}{G(n)\left(1+\frac{1}{\mu(n)^4}\right)}\nonumber\\[5pt]
   &=\frac{G(n-1)}{G(n)}\cdot\frac{1-\frac{1}{\mu(n-1)^4}}{1+\frac{1}{\mu(n)^4}}+\frac{G(n+1)}{G(n)}\cdot\frac{1-\frac{1}{\mu(n+1)^4}}{1+\frac{1}{\mu(n)^4}}.
\end{align}
We proceed to show that for $\mu(n)\geq 5$,
\begin{align}\label{thm-cx-2}
    \frac{1-\frac{1}{\mu(n-1)^4}}{1+\frac{1}{\mu(n)^4}}>1-\frac{6}{\mu(n)^4}
\end{align}
and
\begin{align}\label{thm-cx-3}
    \frac{1-\frac{1}{\mu(n+1)^4}}{1+\frac{1}{\mu(n)^4}}>1-\frac{3}{\mu(n)^4}.
\end{align}
Using \eqref{lem-cx-aaaa}, we derive that for $\mu(n)\geq 3$,
\begin{align}\label{thm-cx-4}
   \frac{1-\frac{1}{\mu(n-1)^4}}{1+\frac{1}{\mu(n)^4}}=\frac{\mu(n)^4\left(\mu(n-1)^4-1\right)}{\left(\mu(n)^4+1\right)\mu(n-1)^4}=\frac{\mu(n)^4\left(\left(\mu(n)^2-\frac{2\pi^2}{3}\right)^2-1\right)}{\left(\mu(n)^4+1\right)\left(\mu(n)^2-\frac{2\pi^2}{3}\right)^2}.
\end{align}
It is easy to show that
\begin{align}\label{thm-cx-5}
   \mu(n)^4\left(\left(\mu(n)^2-\frac{2\pi^2}{3}\right)^2-1\right)&=\mu(n)^8-\frac{4\pi^2}{3}\mu(n)^6+\left(\frac{4\pi^4}{9}-1\right)\mu(n)^4 \nonumber\\[5pt]
   &>\mu(n)^8-\frac{4\pi^2}{3}\mu(n)^6+42\mu(n)^4
\end{align}
and for $\mu(n)\geq 2$,
\begin{align}\label{thm-cx-6}
   \left(\mu(n)^4+1\right)\left(\mu(n)^2-\frac{2\pi^2}{3}\right)^2&=\mu(n)^8-\frac{4\pi^2}{3}\mu(n)^6+\left(\frac{4\pi^4}{9}+1\right)\mu(n)^4-\frac{4\pi^2}{3}\mu(n)^2+\frac{4\pi^4}{9}\nonumber\\[5pt]
   &<\mu(n)^8-\frac{4\pi^2}{3}\mu(n)^6+\left(\frac{4\pi^4}{9}+1\right)\mu(n)^4\nonumber\\[5pt]
   &<\mu(n)^8-\frac{4\pi^2}{3}\mu(n)^6+45\mu(n)^4.
\end{align}
Applying \eqref{thm-cx-5} and \eqref{thm-cx-6} to \eqref{thm-cx-4}, we derive that for $\mu(n)\geq3$,
\begin{align}\label{thm-cx-2aa}
    \frac{1-\frac{1}{\mu(n-1)^4}}{1+\frac{1}{\mu(n)^4}} &>\frac{\mu(n)^8-\frac{4\pi^2}{3}\mu(n)^6+42\mu(n)^4}{\mu(n)^8-\frac{4\pi^2}{3}\mu(n)^6+45\mu(n)^4}\nonumber\\[5pt]
    &=1-\frac{3}{\mu(n)^4-\frac{4\pi^2}{3}\mu(n)^2+45}.
\end{align}
It can be checked that for $\mu(n)\geq 5$,
\[\mu(n)^4-\frac{4\pi^2}{3}\mu(n)^2+45>\frac{1}{2}\mu(n)^4,\]
so we derive from \eqref{thm-cx-2aa} that \eqref{thm-cx-2} holds  for $\mu(n)\geq 5$.
Similarly,
\begin{align}\label{thm-cx-7}
   \frac{1-\frac{1}{\mu(n+1)^4}}{1+\frac{1}{\mu(n)^4}}=\frac{\mu(n)^4\left(\mu(n+1)^4-1\right)}{\left(\mu(n)^4+1\right)\mu(n+1)^4}=\frac{\mu(n)^4\left(\left(\mu(n)^2+\frac{2\pi^2}{3}\right)^2-1\right)}{\left(\mu(n)^4+1\right)\left(\mu(n)^2+\frac{2\pi^2}{3}\right)^2}.
\end{align}
It is evident that
\begin{align}\label{thm-cx-8}
   \mu(n)^4\left(\left(\mu(n)^2+\frac{2\pi^2}{3}\right)^2-1\right)&=\mu(n)^8+\frac{4\pi^2}{3}\mu(n)^6+\left(\frac{4\pi^4}{9}-1\right)\mu(n)^4 \nonumber\\[5pt]
   &>\mu(n)^8+\frac{4\pi^2}{3}\mu(n)^6+42\mu(n)^4
\end{align}
and for $\mu(n)\geq 5$,
\begin{align}\label{thm-cx-9}
   \left(\mu(n)^4+1\right)\left(\mu(n)^2+\frac{2\pi^2}{3}\right)^2&=\mu(n)^8+\frac{4\pi^2}{3}\mu(n)^6+\left(\frac{4\pi^4}{9}+1\right)\mu(n)^4+\frac{4\pi^2}{3}\mu(n)^2+\frac{4\pi^4}{9}\nonumber\\[5pt]
   &<\mu(n)^8+\frac{4\pi^2}{3}\mu(n)^6+45\mu(n)^4.
\end{align}
Substituting \eqref{thm-cx-8} and \eqref{thm-cx-9} to \eqref{thm-cx-7}, we derive that for $\mu(n)\geq5$,
\begin{align*}
    \frac{1-\frac{1}{\mu(n-1)^4}}{1+\frac{1}{\mu(n)^4}} &>\frac{\mu(n)^8+\frac{4\pi^2}{3}\mu(n)^6+42\mu(n)^4}{\mu(n)^8+\frac{4\pi^2}{3}\mu(n)^6+45\mu(n)^4}\\[5pt]
    &=1-\frac{3}{\mu(n)^4+\frac{4\pi^2}{3}\mu(n)^2+45}\\[5pt]
    &>1-\frac{3}{\mu(n)^4}.
\end{align*}
Hence   \eqref{thm-cx-3} holds for  $\mu(n)\geq5$.
Applying \eqref{lem-cx-1}, \eqref{lem-cx-2}, \eqref{thm-cx-2} and \eqref{thm-cx-3} to \eqref{convex-tt}, we derive that   for $\mu(n)\geq 88$,
\begin{align}\label{convex-eqtt}
   &\frac{M_k(n-1)+M_k(n+1)}{M_k(n)} \nonumber\\[5pt]
   &\quad>\left(1+\frac{2\pi^2}{3\mu(n)^2}\right)\left(1-\frac{\pi^2}{\mu(n)^3}\right)\left(1-\frac{\pi^2}{3\mu(n)}+\frac{\pi^4}{18\mu(n)^2}-\frac{18\pi^4+\pi^6}{162\mu(n)^3}\right)\left(1-\frac{6}{\mu(n)^4}\right)\nonumber\\[5pt]
   &\qquad+\left(1-\frac{2\pi^2}{3\mu(n)^2}\right)\left(1+\frac{\pi^2}{3\mu(n)^3}\right)\left(1+\frac{\pi^2}{3\mu(n)}+\frac{\pi^4}{18\mu(n)^2}-\frac{\pi^4}{9\mu(n)^3}\right)\left(1-\frac{3}{\mu(n)^4}\right)\nonumber\\[5pt]
   &\quad >2+\frac{\pi^4}{9\mu(n)^2}-\frac{78}{\mu(n)^3}+\frac{34}{\mu(n)^4}-\frac{152}{\mu(n)^5}+\frac{203}{\mu(n)^6}-\frac{92}{\mu(n)^7}+\frac{988}{\mu(n)^8}+\frac{1169}{\mu(n)^9}\nonumber\\[5pt]
   &\qquad-\frac{1954}{\mu(n)^{10}}+\frac{2459}{\mu(n)^{11}}-\frac{7233}{\mu(n)^{12}}.
\end{align}
It can be checked that for $\mu(n)\geq 8$,
\begin{align*}\left\{
	\begin{aligned}
		&	\frac{\pi^4}{9\mu(n)^2}-\frac{78}{\mu(n)^3}> 0,\\[5pt]
		&	\frac{34}{\mu(n)^4}-\frac{152}{\mu(n)^5}> 0,\\[5pt]
  &\frac{203}{\mu(n)^6}-\frac{92}{\mu(n)^7}>0, \\[5pt]
  &\frac{1169}{\mu(n)^9}-\frac{1954}{\mu(n)^{10}}>0,\\[5pt]
  &\frac{2459}{\mu(n)^{11}}-\frac{7233}{\mu(n)^{12}}>0.
	\end{aligned}\right.
\end{align*}
Hence, it follows from \eqref{convex-eqtt} that   \eqref{convex-eq} holds  for $\mu(n)\geq 88$ (or equivalently, $n\geq 1180$).
It can be checked that \eqref{convex-eq} also holds for $39 \leq n \leq 1180$ if $k=0$ and for $38 \leq n \leq 1180$ if $k=1$. This completes the proof of Theorem \ref{thm-convex}.
\qed

\section{Proof of Theorem \ref{thm-bound-Yk} }

This section is devoted to establishing the upper bound and the lower bound of $M_k(n-1)M_k(n+1)/M_k(n)^2$.

\noindent{\it Proof of Theorem \ref{thm-bound-Yk}. }
Recall that
\[Y_k(n):=\frac{M_k(n-1)M_k(n+1)}{M_k(n)^2} \]
and
\[G(n)=\frac{\pi^2}{12\sqrt{3}\mu(n)^2}\left(1-\frac{1}{\mu(n)}\right)e^{\mu(n)}.\]
Using Theorem \ref{thm-Mk-eq}, we see that for $\mu(n)\geq 88$,
\begin{equation}\label{eq-Yk}
	X(n)L_Y(n)\le Y_k(n)\le X(n)R_Y(n),
\end{equation}
where
\begin{equation}\label{defi-X}
	X(n)=\frac{G(n-1)G(n+1)}{G(n)^2},
\end{equation}
\begin{equation}\label{defi-Z1}
	L_Y(n)=\frac{\left(1-\frac{1}{\mu(n-1)^6}\right)\left(1-\frac{1}{\mu(n+1)^6}\right)}{\left(1+\frac{1}{\mu(n)^6}\right)^2}
\end{equation}
and
\begin{equation}\label{defi-Z2}
	R_Y(n)=\frac{\left(1+\frac{1}{\mu(n-1)^6}\right)\left(1+  \frac{1}{\mu(n+1)^6}\right)}{\left(1-\frac{1}{\mu(n)^6}\right)^2}.
\end{equation}

To prove Theorem \ref{thm-bound-Yk}, we proceed to estimate $X(n)$,  $L_Y(n)$ and $R_Y(n)$ in terms of $\mu(n)$.
We first consider $X(n)$ given by
\begin{align}\label{equation-X}
	X(n)=\frac{\mu(n)^4}{\mu(n-1)^2\mu(n+1)^2}\frac{\left(1-\frac{1}{\mu(n-1)}\right)\left(1-\frac{1}{\mu(n+1)}\right)}{\left(1-\frac{1}{\mu(n)}\right)^2}e^{\mu(n-1)+\mu(n+1)-2\mu(n)}.
\end{align}	
Invoking \eqref{lem-cx-aaaa}, we find that
\begin{align}\label{bound-X-1-aa}
	\frac{\mu(n)^4}{\mu(n-1)^2\mu(n+1)^2}&=\frac{\mu(n)^4}{\left(\mu(n)^2-\frac{2\pi^2}{3}\right)\left(\mu(n)^2+\frac{2\pi^2}{3}\right)}\nonumber\\[5pt]
	&=\frac{\mu(n)^4}{\mu(n)^4-\frac{4\pi^4}{9}}.
\end{align}
It can be calculated that
\[\mu(n)^4-\left(\mu(n)^4-\frac{4\pi^4}{9}\right)\left(1+\frac{4\pi^4}{9\mu(n)^4}+\frac{16\pi^8}{81\mu(n)^8}\right)=\frac{64\pi^{12}}{729\mu(n)^8}>0\]
and for $\mu(n)\ge 8$,
\[\mu(n)^4-\left(\mu(n)^4-\frac{4\pi^4}{9}\right)\left(1+\frac{4\pi^4}{9\mu(n)^4}+\frac{\pi^8}{5\mu(n)^8}\right)=-\frac{\pi^8}{405\mu(n)^4}+\frac{4\pi^{12}}{45\mu(n)^8}<0.\]
Hence, by \eqref{bound-X-1-aa}, we obtain
\begin{equation}\label{bound-X-1}
	1+\frac{4\pi^4}{9\mu(n)^4}+\frac{16\pi^8}{81\mu(n)^8}\leq\frac{\mu(n)^4}{\mu(n-1)^2\mu(n+1)^2}\leq 1+\frac{4\pi^4}{9\mu(n)^4}+\frac{\pi^8}{5\mu(n)^8}.
\end{equation}

We proceed to  estimate the remaining parts on the right-hand side of \eqref{equation-X}. Using \eqref{lem-cx-aaaaccdd},  it is readily checked that   for $\mu(n)\geq 6$,	
\begin{align}
	\tilde{w}(n)<&\mu(n-1)<\hat{w}(n), \label{xn-1}\\[5pt]
	\tilde{y}(n)<&\mu(n+1)<\hat{y}(n), \label{xn+1}
\end{align}
where
\begin{align}
	\begin{aligned}\label{wylabel}
		\tilde{w}(n)&=\mu(n)-\frac{\pi^2}{3\mu(n)}-\frac{\pi^4}{18\mu(n)^3}-\frac{\pi^6}{54\mu(n)^5}-\frac{5\pi^8}{324\mu(n)^7},\\[5pt]
		\hat{w}(n)&=\mu(n)-\frac{\pi^2}{3\mu(n)}-\frac{\pi^4}{18\mu(n)^3}-\frac{\pi^6}{54\mu(n)^5},\\[5pt]
		\tilde{y}(n)&=\mu(n)+\frac{\pi^2}{3\mu(n)}-\frac{\pi^4}{18\mu(n)^3}+\frac{\pi^6}{54\mu(n)^5}-\frac{5\pi^8}{324\mu(n)^7},\\[5pt]
		\hat{y}(n)&=\mu(n)+\frac{\pi^2}{3\mu(n)}-\frac{\pi^4}{18\mu(n)^3}+\frac{\pi^6}{54\mu(n)^5}.
	\end{aligned}
\end{align}
We next show that for $\mu(n)\ge 44$,
\begin{equation}\label{bound-X-2}
	1-\frac{\pi^4}{3\mu(n)^5}-\frac{45}{\mu(n)^6}\leq\frac{\left(1-\frac{1}{\mu(n-1)}\right)\left(1-\frac{1}{\mu(n+1)}\right)}{\left(1-\frac{1}{\mu(n)}\right)^2}\leq 1-\frac{\pi^4}{3\mu(n)^5}.
\end{equation}
Applying \eqref{xn-1} and \eqref{xn+1}, we deduce that for $\mu(n)\geq 6$,
\begin{equation*}
	\frac{\left(1-\frac{1}{\tilde{w}(n)}\right)\left(1-\frac{1}{\tilde{y}(n)}\right)}{\left(1-\frac{1}{\mu(n)}\right)^2}\leq\frac{\left(1-\frac{1}{\mu(n-1)}\right)\left(1-\frac{1}{\mu(n+1)}\right)}{\left(1-\frac{1}{\mu(n)}\right)^2}\leq \frac{\left(1-\frac{1}{\hat{w}(n)}\right)\left(1-\frac{1}{\hat{y}(n)}\right)}{\left(1-\frac{1}{\mu(n)}\right)^2}.
\end{equation*}
Hence to show \eqref{bound-X-2}, it is enough to show that
for $\mu(n)\geq 44$,
\begin{equation}\label{lo-/x}
	\frac{\left(1-\frac{1}{\tilde{w}(n)}\right)\left(1-\frac{1}{\tilde{y}(n)}\right)}{\left(1-\frac{1}{\mu(n)}\right)^2}- \left(1-\frac{\pi^4}{3\mu(n)^5}-\frac{45}{\mu(n)^6}\right)\ge 0
\end{equation}
and
\begin{equation}\label{up-/x}
	\frac{\left(1-\frac{1}{\hat{w}(n)}\right)\left(1-\frac{1}{\hat{y}(n)}\right)}{\left(1-\frac{1}{\mu(n)}\right)^2}-\left( 1-\frac{\pi^4}{3\mu(n)^5}\right)\le 0.
\end{equation}
Observe that
\begin{equation}
	\frac{\left(1-\frac{1}{\tilde{w}(n)}\right)\left(1-\frac{1}{\tilde{y}(n)}\right)}{\left(1-\frac{1}{\mu(n)}\right)^2}- \left(1-\frac{\pi^4}{3\mu(n)^5}-\frac{45}{\mu(n)^6}\right)=\frac{\varphi_1(\mu(n))}{\mu(n)^{20}(\mu(n)-1)^2\tilde{w}(n)\tilde{y}(n)},
\end{equation}
where
\begin{align*}
	\varphi_1(\mu(n))&=\mu(n)^{22}(\tilde{w}(n)-1)(\tilde{y}(n)-1)\\[5pt]
	&\quad	-\mu(n)^{14}\left(\mu(n)^6-\frac{\pi^4}{3}\mu(n)-45\right)(\mu(n)-1)^2\tilde{w}(n)\tilde{y}(n).
\end{align*}
We claim that for $\mu(n)\geq 44$,
\begin{equation}\label{eq-phi}
	\varphi_1(\mu(n))\geq 0.
\end{equation}
It can be checked that $\varphi_1(\mu(n))$ is a polynomial in $\mu(n)$ with degree 18, so we could express
\begin{equation*}	\varphi_1(\mu(n))=\sum_{j=0}^{18}a_j\mu(n)^j.
\end{equation*}
Clearly,
\begin{equation*}
	\varphi_1(\mu(n))\geq -\sum_{j=0}^{16}|a_j|\mu(n)^j+a_{17}\mu(n)^{17}+a_{18}\mu(n)^{18}.
\end{equation*}
Moreover, numerical evidence indicates that for  $0\le j\le 15$ and $\mu(n)\geq 27$,
\begin{equation*}
	-|a_j|\mu(n)^j\geq -|a_{16}|\mu(n)^{16}
\end{equation*}
and
\begin{align*}
	a_{16}=45,\quad
	a_{17}=-90+\frac{\pi^4}{3},\quad
	a_{18}=45-\frac{4\pi^4}{9}.
\end{align*}
It is readily checked that for $\mu(n)\ge 44$,
\[a_{18}\mu(n)^2+a_{17}\mu(n)-17|a_{16}|\geq 0.\]
Assembling all these results above, we conclude that for $\mu(n)\ge 44$,
\begin{equation*}
	\varphi_1(\mu(n))\geq \left(a_{18}\mu(n)^2+a_{17}\mu(n)-17|a_{16}|\right)\mu(n)^{16}\geq 0.
\end{equation*}
This proves \eqref{eq-phi} and so \eqref{lo-/x} is valid. Similarly, observe that
\begin{equation}\label{eq-psi-1}
	\frac{\left(1-\frac{1}{\hat{w}(n)}\right)\left(1-\frac{1}{\hat{y}(n)}\right)}{\left(1-\frac{1}{\mu(n)}\right)^2}-\left( 1-\frac{\pi^4}{3\mu(n)^5}\right)=\frac{-\varphi_2(\mu(n))}{8748\mu(n)^{15}(\mu(n)-1)^2\hat{w}(n)\hat{y}(n)},
\end{equation}
where
\begin{align*}
\varphi_2(s)&=3888\pi^4s^{13}-2916\pi^4s^{12}+810\pi^8s^{10}-1377\pi^8s^9+648\pi^8s^8\\[5pt]
&\quad+33\pi^{12}s^{6}-57\pi^{12}s^5+27\pi^{12}s^4+\pi^{16}s^2-2\pi^{16}s+\pi^{16}.
\end{align*}
It can be readily checked that for $s\geq 2$
\begin{align}\left\{
	\begin{aligned}
	&3888\pi^4s^{13}-2916\pi^4s^{12}\geq 0,\nonumber\\[5pt]
		&810\pi^8s^{10}-1377\pi^8s^9\geq 0,\nonumber\\[5pt]
		&	33\pi^{12}s^{6}-57\pi^{12}s^5\geq 0,\nonumber\\[5pt]
		&\pi^{16}s^2-2\pi^{16}s\geq 0,\nonumber
	\end{aligned}\right.
\end{align}
which implies that  $\varphi_2(s)\geq0$ for $s\geq 2$, thus we have that for $\mu(n)\ge 2$,
\begin{equation}\label{eq-psi-2}
	\varphi_2(\mu(n))\geq 0.
\end{equation}
Hence \eqref{up-/x} is confirmed by applying \eqref{eq-psi-2} to \eqref{eq-psi-1}. Combining \eqref{lo-/x} and \eqref{up-/x}, we obtain \eqref{bound-X-2}.

We proceed to estimate $\exp(\mu(n-1)+\mu(n+1)-2\mu(n))$. Applying \eqref{xn-1}--\eqref{wylabel}, we find that for $\mu(n) \geq 6$,
$$
-\frac{ \pi^4}{9 \mu(n)^3}-\frac{5  \pi^8}{162 \mu(n)^7}<\mu(n-1)+\mu(n+1)-2 \mu(n)<-\frac{ \pi^4}{9 \mu(n)^3}.
$$
It follows that
\begin{align}\label{eq-e1}
    \exp\left(-\frac{ \pi^4}{9 \mu(n)^3}-\frac{5  \pi^8}{162 \mu(n)^7}\right)<\exp\left(\mu(n-1)+\mu(n+1)-2 \mu(n)\right)<\exp\left(-\frac{ \pi^4}{9 \mu(n)^3}\right).
\end{align}
Since for $s<0$,
\[
1+s<e^s<1+s+s^2,\]
we derive that
\begin{align}\label{eq-e2}
\exp\left(-\frac{ \pi^4}{9 \mu(n)^3}\right)<1-\frac{ \pi^4}{9 \mu(n)^3}+\frac{ \pi^8}{81 \mu(n)^6}
\end{align}
and
\begin{align}\label{eq-e3}
\exp\left(-\frac{ \pi^4}{9 \mu(n)^3}-\frac{5  \pi^8}{162 \mu(n)^7}\right)>1-\frac{ \pi^4}{9 \mu(n)^3}-\frac{5  \pi^8}{162 \mu(n)^7}.
\end{align}
Applying \eqref{eq-e2} and \eqref{eq-e3} to  \eqref{eq-e1}, we derive  that for $\mu(n) \geq 6$,
\begin{equation}\label{bound-X-3}
	1-\frac{\pi^4}{9\mu(n)^3}-\frac{5\pi^8}{162\mu(n)^7}< \exp(\mu(n-1)+\mu(n+1)-2\mu(n))< 1-\frac{\pi^4}{9\mu(n)^3}+\frac{\pi^8}{81\mu(n)^6}.
\end{equation}

Applying \eqref{bound-X-1}, \eqref{bound-X-2} and  \eqref{bound-X-3} to  \eqref{equation-X}, we obtain that for $k=0$ or $1$ and $\mu(n)\geq 44 $,
\begin{align}\label{eq-upper-X}
	X(n)\leq\left(1+\frac{4\pi^4}{9\mu(n)^4}+\frac{\pi^8}{5\mu(n)^8}\right)\left(1-\frac{\pi^4}{3\mu(n)^5}\right)\left(1-\frac{\pi^4}{9\mu(n)^3}+\frac{\pi^8}{81\mu(n)^6}\right)
\end{align}	
and	
\begin{align}\label{eq-lower-X}
	X(n)\geq\left(1+\frac{4\pi^4}{9\mu(n)^4}+\frac{16\pi^8}{81\mu(n)^8}\right)\left(1-\frac{\pi^4}{3\mu(n)^5}-\frac{45}{\mu(n)^6}\right)\left(1-\frac{\pi^4}{9\mu(n)^3}-\frac{5\pi^8}{162\mu(n)^7}\right).
\end{align}

 Finally we estimate $L_Y(n)$ and $R_Y(n)$. We claim that for $\mu(n)\ge 16$,
\begin{align}\label{eq-Z}
	L_Y(n)\geq1-\frac{5}{\mu(n)^6} \quad\text {and}\quad  	R_Y(n)\leq 1+\frac{5}{\mu(n)^6}.
\end{align}	
Invoking  \eqref{lem-cx-aaaa}, we  obtain that
$$
L_Y(n)=\frac{\mu(n)^{12}\left(\left(\mu(n)^2+\frac{2 \pi^2}{3}\right)^3-1\right)\left(\left(\mu(n)^2-\frac{2 \pi^2}{3}\right)^3-1\right)}{\left(\mu(n)^6+1\right)^2\left(\mu(n)^4-\frac{4 \pi^4}{9}\right)^3}
$$
and
$$
R_Y(n)=\frac{\mu(n)^{12}\left(\left(\mu(n)^2-\frac{2 \pi^2}{3}\right)^3+1\right)\left(\left(\mu(n)^2+\frac{2 \pi^2}{3}\right)^3+1\right)}{\left(\mu(n)^6-1\right)^2\left(\mu(n)^4-\frac{4 \pi^4}{9}\right)^3} .
$$
Hence,
\begin{equation}\label{eq-ly}
    L_Y(n)-\left(1-\frac{5}{\mu(n)^6}\right)=\frac{\psi_1\left(
\mu(n)\right)}{\mu(n)^6\left(9 \mu(n)^4-4 \pi^4\right)^3\left(\mu(n)^6+1\right)^2}
\end{equation}
and
\begin{equation}\label{eq-ry}
    R_Y(n)-\left(1+\frac{5}{\mu(n)^6}\right)=\frac{-\psi_2\left(\mu(n)\right)}{\mu(n)^6\left(9 \mu(n)^4-4 \pi^4\right)^3\left(\mu(n)^6-1\right)^2},
\end{equation}
where
\begin{align*}
\psi_1(s)&=  729 s^{24}-4860 \pi^4 s^{20}+7290 s^{18}+1296 \pi^8 s^{16}-8748 \pi^4 s^{14}+\left(3645-192 \pi^{12} \right)s^{12} \\[5pt]
&\quad +3888 \pi^8 s^{10}-4860 \pi^4 s^8-576 \pi^{12} s^6+2160 \pi^8 s^4-320 \pi^{12}
\end{align*}
and
\begin{align*}
\psi_2(s)&=729 s^{24}-4860 \pi^4 s^{20}-7290 s^{18}+1296 \pi^8 s^{16}+8748 \pi^4 s^{14}+\left(3645-192 \pi^{12} \right) s^{12} \\[5pt]
&\quad -3888 \pi^8 s^{10}-4860 \pi^4 s^8+576 \pi^{12} s^6+2160 \pi^8 s^4-320 \pi^{12}.
\end{align*}
It can be readily checked that for $s \geq 16$,
\begin{align}\left\{
	\begin{aligned}
&729 s^{24}-4860 \pi^4 s^{20}-7290 s^{18} \geq 0, \nonumber\\[5pt]
&8748 \pi^4 s^{14}+\left(3645-192 \pi^{12} -3888 \pi^8 -4860 \pi^4\right) s^{12} \geq 0, \nonumber\\[5pt]
&2160 \pi^8 s^4-320 \pi^{12}  \geq 0,\nonumber
\end{aligned}\right.
\end{align}
which implies that for $s \geq 16$,
$$
\psi_2(s) \geq 0.
$$
We note that for $s \geq 4$,
$$
\psi_1(s)>\psi_2(s).
$$
Hence for $\mu(n) \geq 16$,
\begin{equation}\label{eq-psi12}
    \psi_1\left(\mu(n)\right) \geq \psi_2\left(\mu(n)\right) \geq 0 .
\end{equation}
and so \eqref{eq-Z} is verified by applying \eqref{eq-psi12} to \eqref{eq-ly} and \eqref{eq-ry} respectively.

Substituting \eqref{eq-upper-X}, \eqref{eq-lower-X} and  \eqref{eq-Z} into \eqref{eq-Yk}, we derive that  for $k=0,\,1$ and $\mu(n)\geq 88$,
\begin{align}\label{bound-Y<}
	Y_k(n)&\leq \left(1+\frac{4\pi^4}{9\mu(n)^4}+\frac{\pi^8}{5\mu(n)^8}\right)\left(1-\frac{\pi^4}{3\mu(n)^5}\right)\nonumber\\[5pt]
	&\quad\times
	\left(1-\frac{\pi^4}{9\mu(n)^3}+\frac{\pi^8}{81\mu(n)^6}\right)\left(1+\frac{5}{\mu(n)^6}\right)
\end{align}	
and	
\begin{align}\label{bound-Y>}
	Y_k(n)&\geq \left(1+\frac{4\pi^4}{9\mu(n)^4}+\frac{16\pi^8}{81\mu(n)^8}\right)\left(1-\frac{\pi^4}{3\mu(n)^5}-\frac{45}{\mu(n)^6}\right)\nonumber\\[5pt]
	&\quad\times \left(1-\frac{\pi^4}{9\mu(n)^3}-\frac{5\pi^8}{162\mu(n)^7}\right)\left(1-\frac{5}{\mu(n)^6}\right).
\end{align}
To prove Theorem \ref{thm-bound-Yk}, it is enough to show that for  $\mu(n)\geq115$,
\begin{align}\label{bound1-Y<}
	&\left(1+\frac{4\pi^4}{9\mu(n)^4}+\frac{\pi^8}{5\mu(n)^8}\right)\left(1-\frac{\pi^4}{3\mu(n)^5}\right)\nonumber\\[5pt]
	&\qquad\times
	\left(1-\frac{\pi^4}{9\mu(n)^3}+\frac{\pi^8}{81\mu(n)^6}\right)\left(1+\frac{5}{\mu(n)^6}\right)\nonumber\\[5pt]
	&\quad<1-\frac{\pi^4}{9\mu(n)^3}+\frac{4\pi^4}{9\mu(n)^4}-\frac{\pi^4}{3\mu(n)^5}+\frac{\frac{\pi^8}{81}+5}{\mu(n)^6}
\end{align}	
and	
\begin{align}\label{bound1-Y>}
	& \left(1+\frac{4\pi^4}{9\mu(n)^4}+\frac{16\pi^8}{81\mu(n)^8}\right)\left(1-\frac{\pi^4}{3\mu(n)^5}-\frac{45}{\mu(n)^6}\right)\nonumber\\[5pt]
	&\qquad\times \left(1-\frac{\pi^4}{9\mu(n)^3}-\frac{5\pi^8}{162\mu(n)^7}\right)\left(1-\frac{5}{\mu(n)^6}\right)\nonumber\\[5pt]
	&\quad>1-\frac{\pi^4}{9\mu(n)^3}+\frac{4\pi^4}{9\mu(n)^4}-\frac{\pi^4}{3\mu(n)^5}-\frac{60}{\mu(n)^6}.
\end{align}
Observe that
\begin{align}
	&\left(1+\frac{4\pi^4}{9\mu(n)^4}+\frac{\pi^8}{5\mu(n)^8}\right)\left(1-\frac{\pi^4}{3\mu(n)^5}\right)\nonumber\\[5pt]
	&\qquad\times
	\left(1-\frac{\pi^4}{9\mu(n)^3}+\frac{\pi^8}{81\mu(n)^6}\right)\left(1+\frac{5}{\mu(n)^6}\right)\nonumber\\[5pt]
	&\qquad-\left(1-\frac{\pi^4}{9\mu(n)^3}+\frac{4\pi^4}{9\mu(n)^4}-\frac{\pi^4}{3\mu(n)^5}+\frac{\frac{\pi^8}{81}+5}{\mu(n)^6}\right)\nonumber\\[5pt]
	&\quad=-\frac{1}{10935\mu(n)^{25}}\sum_{j=0}^{18}b_j\mu(n)^j,\nonumber
\end{align}	
where $b_j$ are real numbers.
Here we just list the values of $b_{16}$, $b_{17}$, $b_{18}$:
\[
b_{16}=6075\pi^4+1620\pi^8,\quad
b_{17}=-2592\pi^8,\quad
b_{18}=540\pi^8.
\]
Clearly
\begin{equation*}
\sum_{j=0}^{18}b_j\mu(n)^j\geq -\sum_{j=0}^{16}|b_j|\mu(n)^j+b_{17}\mu(n)^{17}+b_{18}\mu(n)^{18}.
\end{equation*}
Moreover, it can be checked that for  $0\le j\le 15$ and  $\mu(n)\geq 5$,
\begin{equation*}
	-|b_j|\mu(n)^j\geq -|b_{16}|\mu(n)^{16}
\end{equation*}
and for $\mu(n)\geq 11$,
\[b_{18}\mu(n)^2+b_{17}\mu(n)-17|b_{16}|> 0.\]
Assembling all these results above, we conclude that for $\mu(n)\ge 115$,
\begin{equation*}
	\sum_{j=0}^{18}b_j\mu(n)^j\geq \left(b_{18}\mu(n)^2+b_{17}\mu(n)-17|b_{16}|\right)\mu(n)^{16}> 0.
\end{equation*}	
and so  \eqref{bound1-Y<} is valid.
Similarly, to justify \eqref{bound1-Y>}, we note that
\begin{align}
	& \left(1+\frac{4\pi^4}{9\mu(n)^4}+\frac{16\pi^8}{81\mu(n)^8}\right)\left(1-\frac{\pi^4}{3\mu(n)^5}-\frac{45}{\mu(n)^6}\right)\nonumber\\[5pt]
	&\qquad\times \left(1-\frac{\pi^4}{9\mu(n)^3}-\frac{5\pi^8}{162\mu(n)^7}\right)\left(1-\frac{5}{\mu(n)^6}\right)\nonumber\\[5pt]
	&\qquad-\left(1-\frac{\pi^4}{9\mu(n)^3}+\frac{4\pi^4}{9\mu(n)^4}-\frac{\pi^4}{3\mu(n)^5}-\frac{60}{\mu(n)^6}\right)\nonumber\\[5pt]
	&\quad=\frac{1}{39366\mu(n)^{27}}\sum_{j=0}^{21}c_j\mu(n)^j,\nonumber
\end{align}
where $c_j$ are real numbers.
Here we also list the values of the last three coefficients:
\[
c_{19}=9234\pi^8,\quad
c_{20}=-3159\pi^8,\quad
c_{21}=393660.\]
It is transparent that
\begin{equation*}
\sum_{j=0}^{21}c_j\mu(n)^r\geq -\sum_{j=0}^{19}|c_j|\mu(n)^j+c_{20}\mu(n)^{20}+c_{21}\mu(n)^{21}.
\end{equation*}
Moreover, it can be checked that for  $0\le j\le 18$ and $\mu(n)\geq3$,
\begin{equation*}
	-|c_j|\mu(n)^j\geq -|c_{19}|\mu(n)^{19}
\end{equation*}
and for $\mu(n)\ge 115$,
\[c_{21}\mu(n)^2+c_{20}\mu(n)-20|c_{19}|> 0.\]
Hence we conclude that for $\mu(n)\ge 115$,
\begin{equation*}
	\sum_{j=0}^{21}c_j\mu(n)^j\geq \left(c_{21}\mu(n)^2+c_{20}\mu(n)-20|c_{19}|\right)\mu(n)^{19}> 0,
\end{equation*}	
and so  \eqref{bound1-Y>} is valid.

Substituting  \eqref{bound1-Y<} and \eqref{bound1-Y>} into   \eqref{bound-Y<} and \eqref{bound-Y>}, we arrive at \eqref{thm-bound-Y<} and \eqref{thm-bound-Y>}.  This completes the proof of Theorem \ref{thm-bound-Yk}.
\qed

\section{Proofs of Theorem \ref{thm-LC} and Theorem \ref{thm-highturan} }

In this section, we aim to prove that $M_0(n)$ (resp. $M_1(n)$)  is
log-concave for $n\geq 94$ and satisfies the higher order Tur\'an inequalities for  $n\geq 207$ with the aid of Theorem \ref{thm-bound-Yk}. We first show that  $M_0(n)$ (resp. $M_1(n)$)  is log-concave for $n\geq 94$.

\noindent{\it Proof of Theorem \ref{thm-LC}. }  Recall that
\[Y_k(n):=\frac{M_k(n-1)M_k(n+1)}{M_k(n)^2} \]
To prove Theorem \ref{thm-LC}, it is equivalent to prove that $Y_0(n)<1$ for $n\geq 94$ and $Y_1(n)<1$ for $n\geq 93$.
It is easy to check that for $\mu(n)\geq 4 $,
\[-\frac{\pi^4}{9\mu(n)^3}+\frac{4\pi^4}{9\mu(n)^4}\leq0 \]
and
\[-\frac{\pi^4}{3\mu(n)^5}+\frac{\frac{\pi^8}{81}+5}{\mu(n)^6}\leq 0\]
and by
\eqref{thm-bound-Y<}, we deduce that $Y_k(n)<1$ for $k=0,\,1$  and   $n \geq 2011$. It can be checked that $Y_0(n)<1$ for   $ 94\leq  n \leq2011$  and $Y_1(n)<1$ for   $ 93\leq  n \leq2011$. Hence we conclude that    $ M_0(n)$ is  log-concave for  $n\geq 94$  and  $ M_1(n)$ is  log-concave for  $n\geq 93$.  This completes the proof of Theorem \ref{thm-LC}. \qed

We conclude this paper with the proof of Theorem \ref{thm-highturan} by employing Theorem \ref{thm-bound-Yk}.
The proof of Theorem \ref{thm-highturan} also requires the following lemma given by Jia \cite{Jia-2022}.
\begin{lem}[Jia]\label{lem-Jia}
	Let $u$ and $v$ be two positive real numbers such that $\frac{\sqrt{5}-1}{2}\leq u< v<1$. If
	\[u+\sqrt{(1-u)^3}>v,\]
	then we have
	\[4(1-u)(1-v)-(1-uv)^2>0.\]
\end{lem}

\noindent{\it Proof of Theorem \ref{thm-highturan}. }
To prove $\{M_0(n)\}_{n\geq207}$ and $\{M_1(n)\}_{n\geq206}$  satisfy the  higher order Tur\'an inequalities, it is equivalent to show that
\begin{equation}\label{eq-turan}
	4(1-Y_k(n))(1-Y_k(n+1))-(1-Y_k(n)Y_k(n+1))^2>0
\end{equation}
for $n\geq207 $ if $k=0$  and   for $n\geq206 $ if $k=1$.
We first show that \eqref{eq-turan} holds for $k=0$ or $1$ and $n\geq 2011$.
From Theorem \ref{thm-LC}, we see
that $Y_k(n+1)<1$ for $n\geq 93$. Hence  by Lemma \ref{lem-Jia},   it's enough to show that for $k=0,\,1$ and $n\geq 2011$,
\begin{equation}\label{eq-lem-1}
	\frac{\sqrt{5}-1}{2} \leq Y_k(n)< Y_k(n+1)
\end{equation}
and
\begin{equation}\label{eq-lem-2}
	Y_k(n+1)<Y_k(n)+\sqrt{(1-Y_k(n))^3} .
\end{equation}

Utilizing  \eqref{thm-bound-Y>} in Theorem \ref{thm-bound-Yk} , we see that  for $k=0,\,1$ and $\mu(n)\geq 115$,
\begin{align*}
	Y_k(n)>1-\frac{\pi^4}{9\mu(n)^3}+\frac{4\pi^4}{9\mu(n)^4}-\frac{\pi^4}{3\mu(n)^5}-\frac{60}{\mu(n)^6}.
\end{align*}
Note that for $\mu(n)\geq 4$,
\[\frac{4\pi^4}{9\mu(n)^4}-\frac{\pi^4}{3\mu(n)^5}-\frac{60}{\mu(n)^6}>0\]	
and
\[1-\frac{\pi^4}{9\mu(n)^3}>\frac{\sqrt{5}-1}{2}.\]
Hence we derive that for $\mu(n)\geq 115$,
\[Y_k(n)>\frac{\sqrt{5}-1}{2}.\]
Employing Theorem \ref{thm-bound-Yk} again, we find that  for $k=0,\,1$ and $\mu(n)\geq 115$,
\begin{align}\label{eq-Yn+1}
	Y_k(n+1)-Y_k(n)&>\left(1-\frac{\pi^4}{9\mu(n+1)^3}+\frac{4\pi^4}{9\mu(n+1)^4}-\frac{\pi^4}{3\mu(n+1)^5}-\frac{60}{\mu(n+1)^6}\right)\nonumber\\[5pt]
	&\quad-\left(1-\frac{\pi^4}{9\mu(n)^3}+\frac{4\pi^4}{9\mu(n)^4}-\frac{\pi^4}{3\mu(n)^5}+\frac{\frac{\pi^8}{81}+5}{\mu(n)^6}\right).
\end{align}
Note that for $\mu(n)\geq 3$,
\begin{align}
	\left\{
	\begin{aligned}\label{eq-x}
		&\frac{1}{\mu(n+1)^3}<\frac{1}{\mu(n)^3}-\frac{\pi^2}{2\mu(n)^5},\\[5pt]
		&\frac{1}{\mu(n+1)^4}>\frac{1}{\mu(n)^4}-\frac{4\pi^2}{3\mu(n)^6},\\[5pt]
		&\frac{1}{\mu(n+1)^5}<\frac{1}{\mu(n)^5},\\[5pt]
		&	\frac{1}{\mu(n+1)^6}<\frac{1}{\mu(n)^6}.
	\end{aligned}\right.
\end{align}
Applying \eqref{eq-x} to \eqref{eq-Yn+1}, we attain  that for $k=0,\,1$ and $\mu(n)\geq 115$,
\begin{align}
	Y_k(n+1)-Y_k(n)&>\left(1-\frac{\pi^4}{9\mu(n)^3}+\frac{4\pi^4}{9\mu(n)^4}+\frac{-\frac{\pi^4}{3}+\frac{\pi^6}{18}}{\mu(n)^5}-\frac{\frac{16\pi^6}{27}+60}{\mu(n)^6}\right)\nonumber\\[5pt]
	&\quad-\left(1-\frac{\pi^4}{9\mu(n)^3}+\frac{4\pi^4}{9\mu(n)^4}-\frac{\pi^4}{3\mu(n)^5}+\frac{\frac{\pi^8}{81}+5}{\mu(n)^6}\right)\nonumber\\[5pt]
	&=\frac{\pi^6}{18\mu(n)^5}-\frac{65+\frac{16\pi^6}{27}+\frac{\pi^8}{81}}{\mu(n)^6}.\nonumber
\end{align}
It can be checked that for $\mu(n)\geq 15$,
\[\frac{\pi^6}{18\mu(n)^5}-\frac{65+\frac{16\pi^6}{27}+\frac{\pi^8}{81}}{\mu(n)^6}>0.\]
Hence we derive that for $k=0,\,1$ and $\mu(n)\geq 115$,
\[Y_k(n+1)-Y_k(n)>0,\]
and so \eqref{eq-lem-1} holds for $k=0,\,1$ and $\mu(n)\geq 115$.

To prove  \eqref{eq-lem-2}, invoking  Theorem \ref{thm-bound-Yk} again, we find that for $k=0,\,1$ and $\mu(n)\geq 115$,
\begin{align}\label{eq-Yn+1-Yn-1}
	Y_k(n+1)-Y_k(n)&<\left(1-\frac{\pi^4}{9\mu(n+1)^3}+\frac{4\pi^4}{9\mu(n+1)^4}-\frac{\pi^4}{3\mu(n+1)^5}+\frac{\frac{\pi^8}{81}+5}{\mu(n+1)^6}\right)\nonumber\\[5pt]
	&\quad-\left(1-\frac{\pi^4}{9\mu(n)^3}+\frac{4\pi^4}{9\mu(n)^4}-\frac{\pi^4}{3\mu(n)^5}-\frac{60}{\mu(n)^6}\right)\nonumber\\[5pt]
	&=\frac{\pi^4}{9}\left(\frac{1}{\mu(n)^3}-\frac{1}{\mu(n+1)^3}\right)+\frac{4\pi^4}{9}\left(\frac{1}{\mu(n+1)^4}-\frac{1}{\mu(n)^4}\right)\nonumber\\[5pt]
	&\quad+\frac{\pi^4}{3}\left(\frac{1}{\mu(n)^5}-\frac{1}{\mu(n+1)^5}\right)+\frac{\frac{\pi^8}{81}+5}{x(n+1)^6}+\frac{60}{\mu(n)^6}.
\end{align}
It is easy to check that for $\mu(n)>0$,
\begin{align}
	\left\{
	\begin{aligned}\label{eq-x-2}
		&\frac{1}{\mu(n)^3}-\frac{1}{\mu(n+1)^3}<\frac{\pi^2}{\mu(n)^5}\\[5pt]
		&\frac{1}{\mu(n+1)^4}-\frac{1}{\mu(n)^4}<0\\[5pt]
		&\frac{1}{\mu(n)^5}-\frac{1}{\mu(n+1)^5}<\frac{1}{\mu(n)^5}
	\end{aligned}\right.
\end{align}
and for $\mu(n)\geq 19$,
\begin{align}\label{eq-x-3}
	\frac{\frac{\pi^8}{81}+5}{\mu (n+1)^6}+\frac{60}{\mu(n)^6}<\frac{\frac{\pi^8}{81}+5+60}{\mu(n)^6}<\frac{\pi^2}{\mu(n)^5}.
\end{align}
Applying  \eqref{eq-x-2} and \eqref{eq-x-3} to \eqref{eq-Yn+1-Yn-1}, we deduce  that for $k=0,\,1$ and $\mu(n)\geq 115$,
\begin{align}\label{eq-Yn+1-Yn-2}
	Y_k(n+1)-Y_k(n)&<\frac{\pi^4}{9}\cdot\frac{\pi^2}{\mu(n)^5}+\frac{\pi^4}{3}\cdot\frac{1}{\mu(n)^5}+\frac{\pi^2}{\mu(n)^5}\nonumber\\[5pt]
	&=\frac{\frac{\pi^6}{9}+\frac{\pi^4}{3}+\pi^2}{\mu(n)^5}.
\end{align}
We next show that for $k=0,\,1$ and  $\mu(n)\geq 115$,
\begin{align}\label{eq-Yn+1-Yn-3}
	\sqrt{(1-Y_k(n))^3}>\frac{\frac{\pi^6}{9}+\frac{\pi^4}{3}+\pi^2}{\mu(n)^5}.
\end{align}
From
\eqref{thm-bound-Y<}, we see that  for  $\mu(n)\geq 115$,
\begin{align*}
	1-Y_k(n)&>\frac{\pi^4}{9\mu(n)^3}-\frac{4\pi^4}{9\mu(n)^4}+\frac{\pi^4}{3\mu(n)^5}-\frac{\frac{\pi^8}{81}+5}{\mu(n)^6}\\[5pt]
	&>\frac{\pi^4}{9\mu(n)^3}-\frac{4\pi^4}{9\mu(n)^4}-\frac{\frac{\pi^8}{81}+5}{\mu(n)^4}.
\end{align*}
It can be checked that for $\mu(n)\geq 83$,
\[-\frac{4\pi^4}{9\mu(n)^4}-\frac{\frac{\pi^8}{81}+5}{\mu(n)^4}>-\frac{2}{\mu(n)^3}.\]
It follows that for $\mu(n)\geq 115$,
\begin{align}\label{eq-1-Yk}
	1-Y_k(n)>\frac{\pi^4}{9\mu(n)^3}-\frac{2}{\mu(n)^3}=\frac{\pi^4-18}{9\mu(n)^3}>0.
\end{align}
On the other hand, it is easy to check that for $\mu(n)\geq 33$,
\begin{equation}\label{eq-sqrt-2}
	\frac{\sqrt{(\pi^4-18)^3}}{27\mu(n)^{\frac{9}{2}}}>\frac{\frac{\pi^6}{9}+\frac{\pi^4}{3}+\pi^2}{\mu(n)^5}.
\end{equation}
and so \eqref{eq-Yn+1-Yn-3} holds. In view of Lemma \ref{lem-Jia}, we conclude that \eqref{eq-turan} holds for $k=0$ or $1$ and $n\geq 2011$.   It  can be directly checked that \eqref{eq-turan} is valid  for $207\leq n\leq2010 $ if $k=0$   and   for $206\leq n\leq2010 $ if $k=1$.
This completes the proof of   Theorem \ref{thm-highturan}. \qed

 \vskip 0.2cm
\noindent{\bf Acknowledgment.} This work
was supported by   the National Science Foundation of China.

\vskip 0.2cm

\noindent{\bf Conflict of interest.}  On behalf of all authors, the corresponding author states that there is no conflict of interest. The manuscript has no associated data.


\begin{thebibliography}{0}
	\setlength{\itemsep}{-.8mm}
	\addcontentsline{toc}{section}{}


	\bibitem{Abra-1972}  M. Abramowitz and I.A. Stegun (eds.), Handbook of mathematical functions with formulas, graphs, and mathematical tables, United States Department of Commerce, National Bureau of Standards, 10th printing, 1972.
	
	
	\bibitem{Andrews-Garvan-1988} G.E. Andrews and  F.G.  Garvan, Dyson's crank of a partition, Bull. Amer. Math. Soc. (N.S.) 18(2) (1988) 167--171.

  \bibitem{Andrews-Lewis-2000} G.E. Andrews and R. Lewis, The ranks and cranks of partitions moduli $2$, $3$, and $4$, J. Number Theory, 85 (2000) 74--84.

	
	\bibitem{Bessenrodt-Ono-2016}
	C. Bessenrodt and K. Ono, Maximal multiplicative properties of partitions, Ann. Comb. 20 (1) (2016) 59--64.	

	
	\bibitem{Bringmann-Kane-Rolen-Tripp-2021} K. Bringmann, B. Kane, L. Rolen and Z. Tripp,  Fractional partitions and conjectures of Chern-Fu-Tang and Heim-Neuhauser, Trans. Amer. Math. Soc. Ser. B 8 (2021) 615--634.
	


	
	\bibitem{Chen-2017} W.Y.C. Chen, The spt-function of Andrews, London Math. Soc. Lecture Note Ser., 440
Cambridge University Press, Cambridge, 2017, 141--203.
	
	\bibitem{Chen-Jia-Wang-2019} W.Y.C. Chen, D.X.Q. Jia, and L.X.W. Wang, Higher order
	Tur\'an inequalities for the partition function, Trans. Amer. Math. Soc. 372 (2019) 2143--2165.
	


\bibitem{Choi-Kang-Lovejoy-2009}D. Choi, S.Y. Kang and J. Lovejoy, Partitions weighted by the parity of the crank, J. Combin. Theory Ser. A 116(5) (2009) 1034--1046.


	

	
	
	
	\bibitem{DeSalvo-Pak-2015} S. DeSalvo and I. Pak,  Log-concavity of the partition function, Ramanujan J.  38 (2015) 61--73.
	

	
	\bibitem{Dyson-1944} F.J. Dyson, Some guesses in the theory of partitions, Eureka (1944) 10--15.
	

	
	\bibitem{Garvan-1988}	 F.G. Garvan, New combinatorial interpretations of Ramanujan's partition congruences mod $5$, $7$ and $11$, Trans. Amer. Math. Soc. 305 (1988) 47--77.
	
	\bibitem{Gomez-Zhu-2021}   K. Gomez and E. Zhu, Bounds for coefficients of the $f(q)$ mock theta function and applications to partition ranks, J. Number Theory 226 (2021) 1--23.
	
	\bibitem{Griffin-Ono-Rolen-Zagier-2019} M. Griffin, K. Ono, L. Rolen and D. Zagier, Jensen polynomials for the Riemann zeta function and other sequences, Proc. Natl. Acad. Sci. USA 116 (2019) 11103--11110.
	



 	\bibitem{Gupta-1978} H. Gupta, Finite differences of the partition function, Math. Comp. 32 (1978) 1241--1243.
	



	
\bibitem{Hamakiotes-Kriegman-Tsai-2021} A. Hamakiotes, A. Kriegman and W.-L. Tsai, Asymptotic distribution of the partition crank, Ramanujan J. 56 (2021) 803--820.
	


 \bibitem{Honsberger-1991}  R. Honsberger,
More Mathematical Morsels, Cambridge University Press, 1991.




	
	\bibitem{Jia-2022}
	D.X.Q. Jia, Inequalities for the broken $k$-diamond partition function, J. Number Theory 249 (2023) 314--347.
	



	\bibitem{Lehmer-1939} D.H. Lehmer,  On the remainders and convergence of the series for the partition function, Trans. Amer. Math. Soc. 46 (1939)  362--373.

\bibitem{LocusDawsey-Masri-2019}	M. Locus Dawsey and R. Masri, Effective bounds for the Andrews spt-function, Forum Math. 31(3) (2019) 743--767.


\bibitem{Masri-2021}  J. Males, Asymptotic equidistribution and convexity for partition ranks,
 Ramanujan J. 54 (2021) 397--413.
	
	\bibitem{Masri-2016} R. Masri, Singular moduli and the distribution of partition ranks modulo $2$, Math. Proc. Cambridge Philos. Soc. 160 (2016) 209--232.
	

	
	\bibitem{Ono-Pujahari-Rolen-2019} K. Ono, S. Pujahari and L. Rolen, Tur\'an inequalities for the plane partition function, Adv. Math. 409 (2022) Paper No. 108692, 31 pp.
	
	
	\bibitem{Zapata-Rolon-2013}J.M. Zapata Rol\'on, Asymptotische Werte von Crank-Differenzen (Asymptotic values of crank differences), Diploma thesis (2013).
	

\end{thebibliography}
\end{document}